\documentclass{article}
\usepackage{amsmath,amssymb}

\numberwithin{equation}{section}
\newtheorem{prop}{Proposition}

\newtheorem{lem}{Lemma}
\newtheorem{thm}{Theorem}
\newcommand{\itg}[0]{\mathbb{Z}}              %rational integers
\newcommand{\cpx}[0]{\mathbb{C}}              %Complex number field
\newcommand{\path}[3]{\widehat{\cal P}({#1})^{#2}_{#3}} %space of path
\newcommand{\la}{\lambda}
\newcommand{\ep}[0]{\varepsilon}

\newcommand{\al}[0]{\alpha}
\newcommand{\be}[0]{\beta}
\newcommand{\ka}[0]{\kappa}
\newcommand{\h}[0]{\hbar}
\newcommand{\PI}[0]{2\pi i}
\newcommand{\uu}[0]{\widetilde{u}}
\newcommand{\vv}[0]{\widetilde{v}}
\newcommand{\ww}[0]{\widetilde{w}}
\newenvironment{prf}{\noindent{\it Proof\/}.}{$\;\square$\par\medskip}
\newenvironment{prf2}{\noindent{\it Proof of Theorem 2\/}.}{$\;\square$\par\medskip}

\begin{document}

\title
{Commuting difference operators \\ arising from the elliptic $C_2^{(1)}$ -face model}
\author
{
Koji HASEGAWA${}^{1}$\thanks{{\it E-mail adress} : kojihas@math.tohoku.ac.jp} , 
Takeshi IKEDA${}^{2}$\thanks{{\it E-mail adress} : ike@xmath.ous.ac.jp}
and
Tetsuya KIKUCHI${}^{3}$\thanks{{\it E-mail adress} : tkikuchi@math.tohoku.ac.jp} 
}
\date{}

\maketitle
%\vskip-.5\baselineskip
\centerline{${}^{1,3}${\it Mathematical Institute, Tohoku University,}}
\centerline{\it Sendai 980-8578, JAPAN}
\medskip\centerline{${}^{2}${\it Department of Applied Mathematics, Okayama University of Science, }}
\centerline{\it Okayama 700-0005, JAPAN 
}

\begin{abstract}
\noindent
We study a pair of commuting difference operators 
arising from the elliptic 
$C_2^{(1)}$-face model.
The operators, whose coefficients are expressed 
in terms of the 
Jacobi's elliptic theta function, 
act on the space of meromorphic functions on 
the weight space 
of the $C_2$ type simple Lie algebra.
We show that the space of 
functions spanned by 
the level one characters of the affine Lie algebra 
$\widehat{\mathfrak{sp}}(4,\cpx)$
is invariant under the action of the difference operators.
\end{abstract}

\section{Introduction}
In Ref.\cite{has}, 
one of the authors constructed an L-operator
for Belavin's elliptic quantum R-matrix \cite{Belavin}
acting on the space of meromorphic functions on the weight space 
of the $A_n$ type simple Lie algebra.
The traces of the L-operator, the transfer matrices, give rise to a
family of commuting difference operators
with elliptic theta function coefficient. 
In Ref.\cite{has97}, they are
actually equivalent to Ruijsenaars' operators \cite{Ruij},
which are elliptic extension of Macdonald's $q$-difference operators \cite{Macbook}.
The aim of the present paper is to take 
a step toward a generalization of the above construction 
to the root systems other than the type $A$. 
In this paper, we construct
a pair of commuting difference operators 
acting on the space of functions on the $C_2$ type weight space. 

In the construction of Refs.\cite{has} and \cite{has97},
a relation between Belavin's elliptic quantum R-matrix
and the face-type solution of the Yang-Baxter equation (YBE)
\cite{JMO2}, especially the {\em intertwining vectors} 
\cite{Baxter,JMO88},
played the central role. 
For the root systems
other than type $A$,
it is known
no vertex-type R-matrices
nor the intertwining vectors.
Nevertheless, the face-type solutions of the YBE are known for all classical 
Lie algebras and their vector representations \cite{JMO2}. 
We will utilize this type of 
solution to introduce
the difference operators.
We take
traces (see the section \ref{difference}) of 
the fused Boltzmann weights
to obtain a pair of difference operators (Theorem \ref{commute}).

We also show that the space which is spanned by 
the level one characters of the affine Lie algebra 
$\widehat{\mathfrak{sp}}(4,\cpx)$
is invariant under the action of the difference operators
(Theorem \ref{theta}).

The plan of this paper is as follows.
In section \ref{Results}, we prepare the notation used in the text 
and state the main results.
In section \ref{DefOfFace}, we review the $C_n^{(1)}$-face model 
\cite{JMO2}
in the vector representation, which was given by a set of functions
called Boltzmann weights.
In section \ref{Fusion}, we introduce
the {\em path space}, on which the set of Boltzmann weights 
act naturally as linear maps and thereby explain the notion of
so-called {\em fusion procedure} (see for example Ref.\cite{has97} and 
references there in). 
We also give a set of formula for {\em fused} Boltzmann weights,
which leads to the explicit formula of our difference operators
(Theorem\ref{commute},(ii)).
In section \ref{difference}, we prove
the commutativity of the difference operators.
In section \ref{Theta}, we prove a property 
that the difference operators preserve a 
three dimensional subspace spanned by the level one
characters of the affine Lie algebra $\widehat{\mathfrak{sp}}(4,\cpx)$
\cite{KacPeterson}.
In appendix, we give a formula of a 
similarity transformation of the Boltzmann weights.

Our result can be seen as a type $C$ generalization of
Felder and Varchenko's work \cite{FV}, where they showed
that the Ruijsenaars system of difference operators can 
be recovered from the
dynamical R-matrices, which is nothing but the face-type solution of the YBE.

On the other hand, 
a $BC_n$ generalization of Macdonald polynomial theory 
is studied by Koornwinder \cite{Koorn}.
In Ref.\cite{vD1} van Diejen constructed the corresponding family
of $q$-difference
operators and he studied 
its elliptic extension in Ref.\cite{vD2}.
He succeeded in constructing two elliptic commuting operators,
one is of the 1st order and the other is of the $n$-th order,
so that they give rise to 
an elliptic extension of difference quantum
Calogero-Moser system of type $BC_2$ \cite{vD1} 
in $n=2$. 
It is likely that our operators can be identified with his system
with special choice of parameters.
We hope to report on this issue in the near future.

Extending this work by van Diejen,
Hikami and Komori rescently obtained a general family
of $n$-commuting difference operators with elliptic function coefficients \cite{KH97, KH98}.
Besides the step parameter of difference oprators 
and the modulus of elliptic functions,  
the family contains ten arbitrary parameters.
Their construcion uses 
Shibukawa-Ueno's elliptic $R$-operator \cite{SU}
togather with the elliptic $K$-operators \cite{KH_K-op1, KH_K-op2},
the elliptic solution to the reflecion equation, 
and can be regarded as an elliptic 
generalization of Dunkl type operator approach to those systems,  
which have been extensively used by Cherednik \cite{C} 
(see Ref.\cite{N} for $BC_n$ case). 
It would be interesting if one can find an explicit relationship
between their approach and ours.

\section{Notation and results} \label{Results}
Let $\mathfrak{h}$ be a fixed Cartan subalgebra of the simple
Lie algebra $\mathfrak{g}:=\mathfrak{sp}(4,\cpx)$ and 
denote by $\mathfrak{h}^*$ the dual space of $\mathfrak{h}$. 
We realize the root system $R$ for $(\mathfrak{g},\mathfrak{h})$
as $R:=\{ \pm(\ep_1\pm\ep_2), \pm 2\ep_1, \pm 2\ep_2\}\subset \mathfrak{h}^*.$
A normalized Killing form $(\,,\,)$ is given by
$(\ep_j,\ep_k)=\frac{1}{2}\delta_{jk}.$
We will often identify the space $\mathfrak{h}$ and its
dual $\mathfrak{h}^*$ via the form $(\,,\,)$.
The fundamental weights
are given by $\varpi_1=\ep_1,\varpi_2=\ep_1+\ep_2$.
Let ${\cal{P}}_d$ be the set of weights for the fundamental representation
$L(\varpi_d)$.
We have
\begin{equation}
{\cal{P}}_1=\{\pm\ep_1,\pm\ep_2\},\quad
{\cal{P}}_2=\{\pm(\ep_1 \pm \ep_2), 0 \}.\label{Pd}
\end{equation}
Note that, in these cases, the multiplicity of the weights are all one.

Fix an elliptic modulus 
$\tau$ in the upper half plane $\Im\tau>0$ and
a generic nonzero complex number $\h$.
Let $[u]$ denote the Jacobi theta function with elliptic
nome ${\rm p}:=e^{2\pi i\tau}\,(\Im \tau>0)$ defined by
\begin{equation*}
[u]:=i{\rm p}^{1/8}
\sin\pi u\prod_{m=1}^\infty(1-2{\rm p}^m\cos2\pi u+{\rm p}^{2m})(1-{\rm p}^m).
\end{equation*}
This is an odd function and has the following quasi-periodicity
\begin{equation}
[u+m]=(-1)^m[u],\quad
[u+m\tau]=(-1)^me^{-\pi im^2\tau-2\pi imu}[u]\quad(m\in \mathbb Z).
\label{quasi_Per}
\end{equation}

Let $d,d'$ be $1$ or $2.$
Then the $C_2^{(1)}$ type Boltzmann weights of the type $(d,d')$
are given as follows.
For any square 
$\begin{pmatrix}
        \la     &\mu    \\
        \kappa  &\nu    \\
\end{pmatrix}
\;(\la,\mu,\nu,\kappa\in\mathfrak{h}^*)$ 
of weights,
the Boltzmann weight
$W_{dd'}\!
\left(\left.\begin{array}{ll}
        \la     &       \mu     \\
        \kappa  &       \nu     \\
\end{array}\,\right|u
\right)
$
is given by as a function of the spectral parameter $u\in\cpx.$
See the next section for the explicit formula for $W_{11},$
which are expressed by the Jacobi theta function.

They satisfy the condition
\begin{equation*}
W_{dd'}\!
\left(\left.\begin{array}{ll}
        \la     &       \mu     \\
        \kappa  &       \nu     \\
\end{array}\,\right|u
\right)=0
\quad \mbox{unless} \quad \mu-\la,\nu-\kappa\in 2\h{\cal P}_d,\;
\kappa-\la,\nu-\mu \in 2\h{\cal P}_{d'},
\end{equation*}
and solve the YBE
\begin{align}
&\sum_\eta W_{dd'}
\left(\left.\begin{array}{ll}
        \rho    &       \eta    \\
        \sigma  &       \kappa  \\
\end{array}\,\right|u-v
\right)
W_{dd''}
\left(\left.\begin{array}{ll}
        \la     &       \mu     \\
        \rho    &       \eta    \\
\end{array}\,\right|u-w
\right)
W_{d'd''}
\left(\left.\begin{array}{ll}
        \mu     &       \nu     \\
        \eta    &       \kappa  \\
\end{array}\,\right|v-w
\right)\nonumber\\
=&
\sum_\eta W_{d'd''}
\left(\left.\begin{array}{ll}
        \la     &       \eta    \\
        \rho    &       \sigma  \\
\end{array}\,\right|v-w
\right)
W_{dd''}
\left(\left.\begin{array}{ll}
        \eta    &       \nu     \\
        \sigma  &       \kappa  \\
\end{array}\,\right|u-w
\right)
W_{dd'}
\left(\left.\begin{array}{ll}
        \la     &       \mu     \\
        \eta    &       \nu     \\
\end{array}\,\right|u-v
\right).
\label{YBE1}
\end{align}
The original Boltzmann weights in Ref.\cite{JMO2} are of the type $(1,1)$
in the above terminology.
We generalized it by the fusion procedure (see the section \ref{Fusion})
for the present purpose.

For $\la\in\mathfrak{h}^*$ and $p\in{\cal P}_d\,(d=1,2)$,
we put
\begin{equation*}
\la_{p}:=(\la,p)
  \label{lambda_p}.
\end{equation*} 

\begin{thm}\label{commute}
Let $M_d(u)\;(u\in \cpx,d=1,2)$ be the following 
difference operators acting on the space of functions on
$\mathfrak{h}^*$ 
\begin{equation*}
(M_d(u)f)(\la):=\sum_{p\in {\cal{P}}_d}
{W_{d2}}\left(\left.\begin{array}{ll}    
          \la   &        \la+2\h{p}        \\
          \la   &        \la+2\h{p}       
\end{array}\right|u \right)
\;T_{\widehat{p}}f(\la) ,
\end{equation*}
where
$
T_{\widehat{p}}f(\la):=f(\la+2\h{p}\,).
$

$({\rm i})$
We have 
$M_d(u)M_{d'}(v)=M_{d'}(v)M_d(u)\;(u,v\in\cpx,d,d'=1,2).$

$({\rm ii})$
Let us define the following difference operators
independent of the spectral parameter $u$
\begin{equation*}
\widetilde{M}_1 :=
\sum_{p\in{\cal{P}}_1}\prod_{\substack{q\in{\cal{P}}_1 \\ q\ne \pm p}}
\frac{[\la_{p+q}-\h]}{[\la_{p+q}]}
T_{\widehat{p}},
\end{equation*}
\begin{equation*}
\widetilde{M}_2 :=
\sum_{\substack{p=\pm\ep_1 \\ q=\pm\ep_2}}
\left(\frac{[\la_{p+q}-\h]}{[\la_{p+q}+\h]}
T_{\widehat{p}}T_{\widehat{q}}
+\frac{[2\h]}{[6\h]}
\frac{[2\la_p+2\h]}{[2\la_p]}
\frac{[2\la_q+2\h]}{[2\la_q]}
\frac{[\la_{p+q}-5\h]}{[\la_{p+q}+\h]}
\frac{[\la_{p+q}+2\h]}{[\la_{p+q}]}
\right).
\end{equation*}
Then we have 
$M_1(u) = F(u)\widetilde{M}_1, M_2(u) = G(u) (\widetilde{M}_2 - H(u))$, where
\[ F(u) :=  \frac{[u]\,[u+2\h]^2\,[u+4\h]}{[-3\h]^2\,[\h]^2} , \]
\begin{equation}
 G(u) := \frac{[u-\h]\,[u]^2\,[u+\h]\,[u+2\h]\,[u+3\h]^2\,[u+4\h]}{[-3\h]^4\,[\h]^4},\label{G(u)}
\end{equation}
\[ \mbox{and}\qquad H(u) :=\frac{[u+6\h]\,[u-3\h]\,[2\h]}{[u]\,[u+3\h]\,[6\h]}. \]
\end{thm}

In section \ref{Theta}, we introduce a space of Weyl group invariant 
theta functions, which are preserved by the actions 
of the difference operators.
For $\beta\in {\mathfrak h}^*$, we introduce the following 
operators $S_{\tau\beta},S_{\beta}$ acting on the functions on
$\mathfrak{h}^*$: 
\begin{align*}
&(S_{\tau\beta}f)(\la)
:=\exp\left[\PI \left((\la,\beta)+\tau(\beta,\beta)/2\right)\right]
f(\la+\tau\beta), \\
&(S_{\beta}f)(\la):=f(\la+\beta).
\end{align*}
They satisfy Heisenberg's relations
\begin{equation}
S_\beta S_\gamma =S_\gamma S_\beta,\quad S_{\tau\beta} S_{\tau\gamma} =S_{\tau\gamma} S_{\tau\beta},
\quad S_\gamma S_{\tau\beta}=e^{2\pi i(\gamma,\beta)}S_{\tau\beta}S_\gamma \label{Heise}
\end{equation}
($\gamma,\beta,\in {\mathfrak h}^*$).

Let $Q^\vee,P^\vee$ be the coroot and coweight lattice respectively.
Let $W\subset GL(\mathfrak{h}^*) $ denote the Weyl group for
$(\mathfrak{g},\mathfrak{h})$.
Let $Th^W$ be a space of $W$-invariant theta functions
defined by:
\begin{equation*}
Th^W
:=\left\{ f\;\mbox{is a holomorphic function over}\;\mathfrak{h}^* \left|
\begin{matrix} S_{\tau\al}f=S_\al f=f  & (\forall\al\in Q^\vee) \\
f(w\la)=f(\la) & (\forall w\in W)
\end{matrix} \right. \right\}.
\end{equation*}
It is well-known
that the space is spanned by the 
level one characters of the affine Lie algebra 
$\widehat{\mathfrak{sp}}(4,\cpx)$, and the dimension of this space
is three.

\begin{thm}\label{theta}
We have 
\begin{equation*}
\widetilde{M}_d(Th^W)\subset Th^W\;(d=1,2).
\end{equation*}
\end{thm}
The corresponding facts in the case of $A$ type are proved
in Refs.\cite{has93} and \cite{has97}. 

\section{The $C_n^{(1)}$-face model}\label{DefOfFace}
Fix an integer $n\geq 2.$
We review the definition of the $C_n^{(1)}$-face model given in Ref.\cite{JMO2}.
We realize the root system $R$ of the type $C_n$ as 
$$
R:= \{ \pm(\ep_j \pm \ep_k), \pm 2\ep_l \, | \, 1 \leq j<k \leq n ,1 \leq l\leq n \},
$$
where $\{\ep_j\}_{j=1}^n$ is a basis of a complex vector space
denoted by $\mathfrak{h}^*$ with
a bilinear form $(\,,\,)$
defined by
\begin{equation*}
(\ep_j,\ep_k):=\frac{1}{2}\delta_{jk}.
\end{equation*}
The vector space $\mathfrak{h}^*$ can be identified with the  
dual space of a Cartan subalgebra $\mathfrak{h}$ of 
the simple Lie algebra $\mathfrak{sp}(2n,\cpx).$
The fundamental weights $\varpi_j(1\leq j\leq n)$
are given by $\varpi_j=\ep_1+\ep_2+\dots +\ep_j.$
Let ${\cal P}$ denote the set of weights that belongs to
the vector representation $L(\varpi_1)$ of $\mathfrak{sp}(2n,\cpx).$
We have 
\begin{equation*}
{\cal{P}}=\{\pm\ep_1,\pm\ep_2,\dots,\pm\ep_n\}.
\end{equation*}
Note that the multiplicity of the weights in ${\cal P}$
are all one. 

We shall use the following notation frequently,
\begin{equation*}
\widehat{\cal P}:=2\h {\cal P},\quad \mbox{and}\quad
\widehat{p}:=2\h p\;(p\in {\cal P})
\label{coordinate}.
\end{equation*}

The Boltzmann weights are given by a set of functions of spectral parameter 
$u\in\cpx$ defined for 
any square $\begin{pmatrix}\la&\mu\\ \nu&\kappa
\end{pmatrix}$ of elements of
$\mathfrak{h}^*.$
Let us denote the functions by $W
\left(\left.\begin{array}{ll}
          \la   &        \mu        \\
         \kappa &\nu   \\
\end{array}\;\right|u
\right).$
They satisfy the condition
$$W
\left(\left.\begin{array}{ll}
          \la   &        \mu        \\
         \kappa &\nu   \\
\end{array}\;\right|u
\right)=0
\quad \mbox{unless}\quad \mu-\la,\nu-\mu,\kappa-\la,\nu-\kappa\in \widehat{\cal P}.
$$
For $p,q,r,s\in{\cal P}$ such that $p+q=r+s$, we will write 
\begin{equation*}
\begin{matrix}    
                &p              &               \\
   s\!\!\!\!    &\boxed{u}      &\!\!\!\!\! q   \\
                &r              &               \\
\end{matrix}
=
W\left(\left.\begin{array}{ll}
          \la   &        \la+\widehat{p}        \\
         \la+\widehat{s}&\la+\widehat{p}+\widehat{q}    \\
\end{array}
\;\right|u
\right).
\label{face}
\end{equation*}
They are explicitly given as follows:
\begin{align}
\begin{matrix}    
                &p       &      \\
   p\!\!\!\!    &\boxed{u} &\!\!\!\!\! p\\
                &p      &       \\
\end{matrix}&=
\frac{[c-u]\,[u+\h]}{[c]\,[\h]},\label{2-5a}
\\
\begin{matrix}  
                &p       &      \\
          p \!\!\!\!    &\boxed{u} &\!\!\!\!\! q\\
                &q      &       \\
\end{matrix}&=
\frac{[c-u]\,[\la_{p-q}-u]}{[c]\,[\la_{p-q}]} \qquad (p\ne \pm q),\label{2-5b}
\\
\begin{matrix}    
                &q       &      \\
          p \!\!\!\!     &\boxed{u} &\!\!\!\!\! p\\
                &q      &       \\
\end{matrix}&=
\frac{[c-u]\,[u]\,[\la_{p-q}+\h]}{[c]\,[\h]\,[\la_{p-q}]} 
\qquad (p\ne \pm q),\label{2-5c}
\\
\begin{matrix}    
                &q       &      \\
          p \!\!\!\!    &\boxed{u} &\!\!\!\!\! -q\\
                &-p &       \\
\end{matrix}&=
-\frac{[u]\,[\la_{p+q}+\h+c-u]}{[c]\,[\la_{p+q}+\h]} 
\frac{[2\la_p+2\h]}{[2\la_q]}
\frac{\prod_{r\ne \pm p}[\la_{p+r}+\h]}
{\prod_{r\ne \pm q}[\la_{q+r}]} \qquad (p\ne q),\label{2-5d}
\\
\begin{matrix}    
                &p       &      \\
          p \!\!\!\!    &\boxed{u} &\!\!\!\!\! -p\\
                &-p &       \\
\end{matrix}&=
\frac{[c-u]\,[2\la_p+\h-u]}{[c]\,[2\la_p+\h]}
-\frac{[u]\,[2\la_p+\h+c-u]}{[c]\,[2\la_p+\h]}
\frac{[2\la_p+2\h]}{[2\la_p]}
\prod_{q\ne\pm p}
\frac{[\la_{p+q}+\h]}{[\la_{p+q}]}.\label{2-5e}
\end{align}

The {\it crossing parameter} $c$ in the above formulas 
are fixed to be
\begin{equation}
c:=-(n+1)\h.\label{CrossingParam}
\end{equation}

\begin{prop}
The Boltzmann weights (\ref{2-5a},\ref{2-5b},\ref{2-5c},\ref{2-5d},\ref{2-5e})
enjoy the following properties. 
%\begin{align}
%\mbox{\textit{Initial condition:}} & \\

\noindent
Initial condition:
\begin{equation}
\sum_\eta W
\left(\left.\begin{array}{ll}
         \la    &       \mu       \\
         \ka    &       \nu        \\
\end{array}\;\right|0
\right) \;
=\;
\delta_{\mu \kappa}.
\label{Initial}
\end{equation}
Inversion relation:
\begin{equation}
\sum_\eta W
\left(\left.\begin{array}{ll}
         \la    &       \eta       \\
         \ka    &       \nu        \\
\end{array}\;\right|u
\right)
W
\left(\left.\begin{array}{ll}
         \la    &       \mu       \\
         \eta    &       \nu        \\
\end{array}\;\right|-u
\right) \,
=\;
\delta_{\mu \kappa}\;
\frac{[c+u]\,[c-u]\,[\h +u]\,[\h -u]}{[c]^2\,[\h]^2}.
\label{FirstInversion}
\end{equation}
Crossing symmetry:
\begin{equation}
W
\left(\left.\begin{array}{ll}
          \la   &        \mu        \\
         \kappa &\nu   \\
\end{array}\;\right|u
\right) \;
=\;
\frac{g(\la,\ka)}{g(\mu,\nu)}\;
W
\left(\left.\begin{array}{ll}
         \ka    &       \la        \\
         \nu    &       \mu        \\
\end{array}\;\right|c-u
\right),
\label{Crossing}
\end{equation}
where we put
\begin{equation*}
g(\la,\mu):=[2\mu_p]\prod_{\substack{q \in {\cal P} \\ q\ne \pm p}}[\mu_{p+q}]
\quad (\mu=\la+\widehat{p}, \; p\in{\cal P}).
\end{equation*}
Reflection symmetry:
\begin{equation}
W
\left(\left.\begin{array}{ll}
         \la    &       \mu        \\
         \ka    &       \nu        \\
\end{array}\;\right|u
\right) \;
=\;
\frac{g(\la,\ka)g(\ka,\nu)}{g(\la,\mu)g(\mu,\nu)}\;
W
\left(\left.\begin{array}{ll}
         \la    &       \ka        \\
         \mu    &       \nu        \\
\end{array}\;\right|u
\right).
\label{Reflection}
\end{equation}
\end{prop}
\begin{prf}
The equation (\ref{Initial}) is trivial.
The two types of symmetries 
(\ref{Crossing}),(\ref{Reflection}) are easily
checked by the explicit form.
In the case of $\la=\nu$ the equation 
(\ref{FirstInversion}) is reduced to the following
\begin{align}
\sum_{r\in{\cal P}}\,&\frac{[\la_p+\la_r+\h+c-u][\la_q+\la_r+\h+c+u]}{[\la_p+\la_r+\h][\la_q+\la_r+\h]}\,G_{\la r}\nonumber\\
=&\delta_{p,q}\frac{[c-u][c+u][2\la_p][2\la_q+2\h]}{[\h]^2[2\la_p+\h]^2}\,G_{\la p}^{-1}\nonumber\\
+&\frac{[c+u][2\la_p+\h+u][\la_p+\la_q+\h+c-u]}
{[u][2\la_p+\h][\la_p+\la_q+\h]} \nonumber \\
-&\frac{[c-u][2\la_q+\h-u][\la_p+\la_q+\h+c+u]}
{[u][2\la_q+\h][\la_p+\la_q+\h]}.\label{iii}
\end{align}
Here we denote by $G_{\la p}$ the following function
\begin{equation}
G_{\la p} :=-\frac{[2\la_p+2\h]}{[2\la_p]} \prod_{\substack{ r \in {\cal P} \\ r \neq \pm p}} 
\frac{[\la_{p+r}+\h]}{[\la_{p+r}]}\quad (p\in {\cal P}).\label{G}
\end{equation}
One can find a proof of the equation (\ref{iii}) 
in Ref.\cite{JMO2}(see (3.5) and Lemma 3).
The cases $\mu=\la+2\widehat{p}\;(p\in{\cal P})$ are trivial.
The remaining cases are easily checked by using the following 
{\em three term identity}:
\begin{eqnarray}
&[u+x][u-x][v+y][v-y]-[u+y][u-y][v+x][v-x] \nonumber \\
&=[x+y][x-y][u+v][u-v]
\label{3term}
\end{eqnarray}
$(u,v,x,y \in\cpx)$.
\end{prf}
We adopted a slightly different formulas (\ref{2-5c}),(\ref{2-5d}) 
from the original ones (see (\ref{2-5oc}),(\ref{2-5od})) in Ref.\cite{JMO2}.
In Appendix, we will give a similarity transformation 
(\ref{function_s}),(\ref{Rel_to_JMO}) which
transforms our Boltzmann weights into the original ones. 
Thus one has a way to prove the YBE for our Boltzmann weights,
since such a transformation does not destroy the varidity of the YBE.
If we follow this track, however, we must specify the arguments of the square roots
contained in the expressions of the original formulas and the transformation.
This way of proof may require a rather complicated discussion.
In this paper, we will give a proof of the YBE for our Boltzmann weights directly 
without using the similarity transformation.
In fact, our proof here goes quite parallel to the proof
given in Ref.\cite{JMO2}.
\begin{thm}%{\rm \cite{JMO2}}
The Boltzmann weights 
$W
\left(\left.\begin{array}{ll}
          \la   &        \mu        \\
         \kappa &\nu   \\
\end{array}\;\right|u
\right)$
(\ref{2-5a},\ref{2-5b},\ref{2-5c},\ref{2-5d},\ref{2-5e})
solve the YBE (\ref{YBE1})
for $d=d'=d''=1.$
\end{thm}
\begin{prf}
Set
\begin{equation}
X(\la,\mu,\nu,\ka,\sigma,\rho\,|\,
u,v)
:=\sum_\eta W
\left(\left.\begin{array}{ll}
        \rho    &       \eta    \\
        \sigma  &       \kappa  \\
\end{array}\,\right|u
\right)
W
\left(\left.\begin{array}{ll}
        \la     &       \mu     \\
        \rho    &       \eta    \\
\end{array}\,\right|u+v
\right)
W
\left(\left.\begin{array}{ll}
        \mu     &       \nu     \\
        \eta    &       \kappa  \\
\end{array}\,\right|v
\right),
\end{equation}
\begin{equation}
Y(\la,\mu,\nu,\ka,\sigma,\rho\,|\,
u,v)
:=
\sum_\eta W
\left(\left.\begin{array}{ll}
        \la     &       \eta    \\
        \rho    &       \sigma  \\
\end{array}\,\right|v
\right)
W
\left(\left.\begin{array}{ll}
        \eta    &       \nu     \\
        \sigma  &       \kappa  \\
\end{array}\,\right|u+v
\right)
W
\left(\left.\begin{array}{ll}
        \la     &       \mu     \\
        \eta    &       \nu     \\
\end{array}\,\right|u
\right),
\end{equation}
and
%\mbox{and}\quad 
\begin{equation}
Z(\la,\mu,\nu,\kappa,\sigma,\rho\,|\,
u,v):=X(\la,\mu,\nu,\kappa,\sigma,\rho\,|\,
u,v)-Y(\la,\mu,\nu,\kappa,\sigma,\rho\,|\,
u,v).
\end{equation}
Regarding $Z(\la,\mu,\nu,\kappa,\sigma,\rho\,|\,
u,v)$ as a function of $u$, we denote it by $Z(u)$.

The equations (\ref{Initial}) and (\ref{FirstInversion}) implies 
$Z(0)=Z(-v)=0.$
Since we have 
\begin{equation}
Z(\la,\mu,\nu,\ka,\sigma,\rho\,|\,
u,v)=
-\frac{g(\la,\rho)}{g(\nu,\ka)}
\,
Z(\rho,\la,\mu,\nu,\ka,\sigma\,|\,
c-u-v,u) \label{Crossing_Z}
\end{equation}
by (\ref{Crossing}), 
this shows $Z(c-v)=Z(c)=0$ also.
Thus we have found the four zeros at $u=0,-v,c,c-v$ of $Z(u).$
By the exactly same argument in Ref.\cite{JMO2} using the 
quasi-periodicity property of $Z(u)$, 
(\ref{Crossing_Z}) and the following symmetry (this follows from (\ref{Reflection}))
\begin{equation*}
Z(\la,\mu,\nu,\ka,\sigma,\rho\,|\,
u,v)
=
\frac{g(\la,\rho)g(\rho,\sigma)g(\sigma,\ka)}{g(\la,\mu)g(\mu,\nu)g(\nu,\ka)}
\,
Z(\la,\rho,\sigma,\ka,\nu,\mu\,|\,v,u),
\label{Reflection_Z}
\end{equation*}
we can reduce the proof of the YBE to the following 
two special cases:
%\begin{align}
\begin{equation}
Z(\la,\la+\widehat{p},\la+\widehat{p}+\widehat{q},\la+\widehat{p}+\widehat{q}+\widehat{r},
\la+\widehat{q}+\widehat{r},\la+\widehat{r}\,|\,u,v)=0,
\label{Exception1}
\end{equation}
%\quad(r\ne\pm p,\pm q,p\ne\pm q),\label{Exception1}\\ 
where $r\ne\pm p,\pm q,p\ne\pm q$ and
\begin{equation}
Z(\la,\la+\widehat{p},\la,\la+\widehat{p},
\la,\la+\widehat{p}\,|\,u,v)=0.
\label{Exception2}
\end{equation}
%\end{align}

In the case of the equation (\ref{Exception1}),
each side of the YBE contains
only one term, and they are manifestly the same.
A proof of the last case (\ref{Exception2})
can be found in the original literature \cite{JMO2}.
However, since the proof is brief and seems to contain
some typographical errors, we will describe 
details of it in the following for readers' convenience.

We will prove 
$Z(\la,\la+{\widehat p},\la,\la+{\widehat p},\la,\la+{\widehat p}\,|u,v)=0.$
Regarding $Y(\la,\la+{\widehat p},\la,\la+{\widehat p},\la,\la+{\widehat p}\,|u,v)$
as a function of $\la_p$ we denote 
it by $f(\la_p).$ It reads as
\begin{align}
f(\la_p)=&
G_{\la\,p}\frac{[u][v][w]}{[c]^3}
\sum_{q\in {\cal P}}
\frac{[\la_q+\la_p+\h+\uu]
[\la_q+\la_p+\h+\vv]
[\la_q+\la_p+\h+\ww]}
{[\la_q+\la_p+\h]^3}\,
G_{\la\,q}\nonumber
\\
&+G_{\la\,p}^{-1}\,
\frac{[\uu][\vv][\ww]}{[c]^3}
\frac{[2\la_p+\h-u][2\la_p+\h-v][2\la_p+\h-w]}{[2\la_p+\h]^3}
\nonumber\\
&+\sum_{\mbox{cyclic}}\frac{[u][\vv][\ww]}{[c]^3}
\frac{[2\la_p+\h+\uu][2\la_p+\h-v][2\la_p+\h-w]}{[2\la_p+\h]^3}
\nonumber\\
&+G_{\la\,p}\sum_{\mbox{cyclic}}\frac{[\uu][v][w]}{[c]^3}
\frac{[2\la_p+\h-u][2\la_p+\h+\vv][2\la_p+\h+\ww]}{[2\la_p+\h]^3},
\nonumber
\end{align}
where we put $w=c-u-v,\uu=c-u,\vv=c-v,\ww=c-w$ and 
the summation $\sum_{\mbox{cyclic}}$ is over 
the cyclic permutations of the three variables $(u,v,w).$
From the explicit form, one can see that
$X(\la,\la+{\widehat p},\la,\la+{\widehat p},\la,\la+{\widehat p}\,|u,v)
=f(-\la_p-\h).$
We will prove
$f(\la_p)=f(-\la_p-\h).$

Now consider a function
\begin{equation} 
\Phi(z)
:= 
\frac{[z+\la_p+\h+\uu][z+\la_p+\h+\vv][z+\la_p+\h+\ww]}
{[z+\la_p+\h]^3}
\frac{[0]'}{[\h]}
\frac{[2z+2\h]}{[2z+\h]}
\prod_{q\in {\cal P}}
\frac{[z+\la_q+\h]}{[z+\la_q]}.\nonumber
\end{equation}
One sees that $\Phi(z)$ is a
doubly periodic function of the periods $1$ and $\tau$.
Its poles are located at
$z=-\la_p-\h,\la_q(q\in{\cal P}),-\frac{\h}{2}+\omega
(\omega=0,\frac{1}{2},\frac{\tau}{2},\frac{1+\tau}{2}).$
The pole at $z=-\la_p-\h$ is of the second order, and the others
are simple. 

Let $f_i(\la_p)\;(i=1,2,3,4)$ denote the $i$-th 
term of the above function $f(\la_p).$
Since we have 
\begin{equation}
\underset{z=\la_q}{\rm Res}\Phi(z)dz
=-\frac{[\la_q+\la_p+\h+\uu][\la_q+\la_p+\h+\vv][\la_q+\la_p+\h+\ww]}{[\la_q+\la_p+\h]^3}G_{\la\,q}\nonumber
\end{equation}
the relation
$\sum{\rm Res}\Phi(z)dz=0$ implies 
$f_1(\la_p)=a(\la_p)+b(\la_p)$
where we set
\begin{align}
a(\la_p)&:=
G_{\la\, p}\frac{[u][v][w]}{[c]^3}
\sum_\omega\underset{z=-\frac{\h}{2}+\omega}{\rm Res}\Phi(z)dz,\label{def_a}\\
b(\la_p)&:=G_{\la\, p}\frac{[u][v][w]}{[c]^3}
\underset{z=-\la_p-\h}{\rm Res}\Phi(z)dz.\label{def_b}
\end{align}
Here the summation $\sum_\omega$ is over the half periods
$\omega=0,\frac{1}{2},\frac{\tau}{2},\frac{1+\tau}{2}.$

From (\ref{quasi_Per}) and (\ref{CrossingParam}),
we have for $\omega=0,\frac{1}{2},\frac{\tau}{2},\frac{1+\tau}{2}$
\begin{equation}
\underset{z=-\frac{\h}{2}+\omega}{\rm Res}\Phi(z)dz
=\frac{1}{2}
\frac{[\la_p+\frac{\h}{2}+\omega+\uu]
[\la_p+\frac{\h}{2}+\omega+\vv]
[\la_p+\frac{\h}{2}+\omega+\ww]}
{[\la_p+\frac{\h}{2}+\omega]^3}
e^{2\pi i\xi(\omega)},\label{omega}
\end{equation}
where we put
$\xi(0)=\xi(\frac{1}{2})=0,\;
\xi(\frac{\tau}{2})=\xi(\frac{1+\tau}{2})=c.$
Combining (\ref{def_a}), (\ref{omega}) and Lemma 3 in Ref.\cite{JMO2}, we can
verify
\begin{equation}
a(\la_p)+f_4(\la_p)-f_2(-\la_p-\h)=
-a(-\la_p-\h)+f_2(\la_p)-f_4(-\la_p-\h)=0.\label{Last1}
\end{equation}

Set $\phi(u)=\frac{d}{du}\log [u],$ then
the residue $\underset{z=-\la_p-\h}{\rm Res}\Phi(z)dz$
can be expressed as
\begin{align}
\underset{z=-\la_p-\h}{\rm Res}\Phi(z)dz
=G_{\la\,p}^{-1}
&\frac{[\uu][\vv][\ww]}{[0]'[\h]^2}
\frac{[2\la_p+2\h][2\la_p]}{[2\la_p+\h]^2} 
\Biggl( 
\sum_{\mbox{cyclic}}
\phi(\uu)
-3\phi(2\la_p)+3\phi(2\la_p+\h)  %\right. 
\nonumber \\
& {}+\phi(\h)+\sum_{\substack{q\in{\cal{P}} \\ q\ne \pm p}}
\left\{
\phi(-\la_p+\la_q)-\phi(-\la_p+\la_q-\h) 
\right\}
\Biggr).\label{res3}
\end{align}
Since
$\phi(u)$ is an odd function, we have from (\ref{def_b}) and (\ref{res3})
\begin{align}
b(\la_p)-b(-\la_p-\h)
= & -3\frac{[u][v][w]}{[c]^3}
\frac{[\uu][\vv][\ww]}{[0]'[\h]^2}
\frac{[2\la_p+2\h][2\la_p]}{[2\la_p+\h]^2}  \nonumber \\
& \times \left\{
\phi(2\la_p)+\phi(2\la_p+2\h)-2\phi(2\la_p+\h)
\right\}.\label{b}
\end{align}

On the other hand, using the identity (see (\ref{3term}))
\begin{align*}
[2\la_p+\h+\uu]&[2\la_p+\h-v][2\la_p+\h-w]
-[2\la_p+\h-\uu][2\la_p+\h+v][2\la_p+\h+w]  \\
&=[\uu][v][w]\frac{[4\la_p+2\h]}{[2\la_p+\h]}
\end{align*}
and its cyclic permutations of $(u,v,w)$,
we have
\begin{equation}
f_3(\la_p)-f_3(-\la_p-\h)
=3\frac{[u][v][w][\uu][\vv][\ww]}{[c]^3}
\frac{[4\la_p+2\h]}{[2\la_p+\h]^4}.
\label{f3}
\end{equation}
Now from (\ref{b}) and (\ref{f3}), we have 
\begin{equation}
b(\la_p)+f_3(\la_p)=b(-\la_p-\h)+f_3(-\la_p-\h),\label{Last2}
\end{equation}
where we used the following identity (Lemma 4 in \cite{JMO2})
\begin{equation}
\phi(u+\h)+\phi(u-\h)-2\phi(u)
=\frac{[\h]^2[2u][0]'}{[u]^2[u-\h][u+\h]}.
\nonumber
\end{equation}
Combining (\ref{Last1}) and (\ref{Last2}) 
we obtained $f(\la_p)=f(-\la_p-\h).$
\end{prf}
%%%%%%%%%%%%%%% Secton 2 %%%%%%%%%%%%%%%%%

\section{Path space and fusion procedure}\label{Fusion} 
In the previous section we introduced 
the Boltzmann weights $W(u)$ of the type $(1,1)$
and proved that they satisfy the YBE.
In what follows, we treat only the case of $n=2.$
To construct commuting difference operators, we need
the general types
of the Boltzmann weights $W_{dd'}(u)$,
which we call the fused Boltzmann weights.

First let us introduce the notion of the path space.
Let $d=1,2.$
For any $u\in \cpx$ and $\la,\mu\in{\mathfrak h}^*$ such that
$\mu-\la\in 2\h {\cal P}_d$, 
we introduce a formal symbol
\begin{equation*}
g_{\la}^\mu(u):=\left\{\begin{array}{cl}
       e^\mu_\la(u)  &
                  :\;d=1
                  \\
       f^\mu_\la(u) & 
                  :\;d=2
       \\
\end{array}\right..
\end{equation*}
See (\ref{Pd}) for the notation ${\cal P}_1$
and ${\cal P}_2.$
We define the complex vector space
$$
\path     {\varpi_d^u}
        {\mu}
        {\lambda}
:=
\left\{\begin{array}{cl}
        {\cpx}\,{g}^\mu_\la(u)  &
                  : \mu-\la\in 2\h{\cal P}_d
                  \\
        0       & :{\rm otherwise}\\
\end{array}\right.
$$
for each $u\in\cpx$,
and the {\em space of paths}
from $\la$ to $\mu$ of the type $(d_1,\ldots,d_k;u_1,\ldots,u_k)$
%\begin{align}
\begin{equation}
\path     {\varpi^{u_1}_{d_1}\otimes\cdots\otimes\varpi^{u_k}_{d_k}}
        {\nu}
        {\lambda} %\nonumber \\
:=
\bigoplus_{\mu_1,\cdots,\mu_{k-1}\in{\mathfrak{h}}^*}
\path     {\varpi^{u_1}_{d_1}}
        {\mu_1}
        {\lambda}
\otimes
\path     {\varpi^{u_2}_{d_2}}
        {\mu_2}
        {\mu_1}
\otimes \cdots \otimes
%\path     {\varpi^{u_{k-1}}_{d_{k-1}}} 
%        {\mu_{k-1}}
%        {\mu_{k-2}}
%        \otimes
\path     {\varpi^{u_{k}}_{d_k}}
        {\nu}
        {\mu_{k-1}}.\label{SpaceOfPath}
%\end{align}
\end{equation}
The following set 
$$\{g_\la^{\mu_1}(u_1)\otimes g_{\mu_1}^{\mu_2}(u_2)
\otimes\cdots 
%\otimes g_{\mu_{k-2}}^{\mu_{k-1}}(u_{k-1})
\otimes g_{\mu_{k-1}}^{\nu}(u_k)\;|\;
\mu_i-\mu_{i-1}\in 2\h{\cal P}_{d_i}(1\leq i\leq k),\mu_0=\la,\mu_k=\nu
\}$$ of {\em paths} 
forms a basis of the space (\ref{SpaceOfPath}).
Set also
$$
\path     {\varpi^{u_1}_{d_1}\otimes\cdots\otimes\varpi^{u_k}_{d_k}}
        {}
        {\lambda}
:=
\bigoplus_
        {\nu\in{\mathfrak h}^*}
\path    {\varpi^{u_1}_{d_1}\otimes\cdots\otimes\varpi^{u_k}_{d_k}}
        {\nu}
        {\lambda}
$$
and
$$
\path   {\varpi^{u_1}_{d_1}\otimes\cdots\otimes\varpi^{u_k}_{d_k}}
        {}
        {}
:=
\bigoplus_
        {\la\in{\mathfrak h}^*}
        \path    {\varpi^{u_1}_{d_1}\otimes\cdots\otimes\varpi^{u_k}_{d_k}}
        {}
        {\lambda}.
$$

In the following, we will construct
the linear operators
$$
W_{dd'}(u-v):\widehat{\cal P}(\varpi_d^u\otimes \varpi_{d'}^v)
\rightarrow
\widehat{\cal P}( \varpi_{d'}^v\otimes\varpi_d^u)
$$
which satisfy the following YBE $(d,d',d''=1,2)$
\begin{eqnarray}
&&\left({\rm id}\otimes W_{dd'}(u-v) \right)
\left(W_{dd''}(u-w)\otimes {\rm id}\right)
\left({\rm id} \otimes W_{d'd''}(v-w) \right) \nonumber\\
&=&\left(W_{d'd''}(v-w)\otimes {\rm id} \right)
\left({\rm id} \otimes W_{dd''}(u-w) \right)
\left(W_{dd'}(u-v)\otimes {\rm id} \right) \nonumber\\
&&:\widehat{\cal P}(\varpi_{d}^u\otimes\varpi_{d'}^v\otimes\varpi_{d''}^w)
\rightarrow
\widehat{\cal P}(\varpi_{d''}^w\otimes\varpi_{d'}^v\otimes\varpi_{d}^u). 
\label{YBE_op}
\end{eqnarray}
First we define a linear operator
$W(\varpi_1^u,\varpi_1^v):
\widehat{\cal P}(\varpi_1^u\otimes \varpi_1^v)
\rightarrow\widehat{\cal P}(\varpi_1^v\otimes \varpi_1^u)$
by 
$$W(\varpi_1^u,\varpi_1^v)\,
e^\mu_\la(u)\otimes e^\nu_\mu(v)
:=\sum_{\kappa\in \mathfrak{h}^*}
W
\left(\left.\begin{array}{ll}
          \la   &        \mu        \\
         \kappa &\nu   \\
\end{array}\;\right|u-v
\right)\,e^\kappa_\la(v)\otimes e^\nu_\kappa(u).
$$
Put $W_{11}(u-v):=W(\varpi_1^u,\varpi_1^v)$,
then the YBE (\ref{YBE_op}) for $d=d'=d''=1$ is 
nothing but (\ref{YBE1}).

To construct $W_{dd'}(u-v)$ other than $W_{11}(u-v),$
we will formulate the fusion procedure.
Put 
\begin{align*}
&W\left(\varpi_1^{u_1}\otimes \varpi_1^{u_2}\otimes \cdots
\otimes \varpi_1^{u_k},\varpi_1^v
\right)  \\
:=&W^{1,2}\left(\varpi_1^{u_1},\varpi_1^v
\right)
W^{2,3}\left(\varpi_1^{u_2},\varpi_1^v
\right)
\cdots
W^{k,k+1}\left(\varpi_1^{u_k},\varpi_1^v
\right)  \label{hfusion} \\
: \widehat{\cal P}(\varpi_1^{u_1}&\otimes\varpi_1^{u_2}\otimes
\dots\otimes\varpi_1^{u_k}\otimes\varpi_1^{v})
\rightarrow
\widehat{\cal P}(\varpi_1^{v}\otimes\varpi_1^{u_1}\otimes\varpi_1^{u_2}
\dots\otimes\varpi_1^{u_k}),\nonumber
\end{align*}
where
\begin{align*}
W&\left(\varpi_1^{u_j},\varpi_1^v
\right)^{j,j+1}
:={\rm id}^{\otimes(j-1)}
\otimes
W\left(\varpi_1^{u_j},\varpi_1^v
\right)
\otimes
{\rm id}^{\otimes(k-j)}\\
:
\widehat{\cal P}&
(\varpi_1^{u_1}\otimes\cdots\otimes\varpi_1^{u_{j-1}}\otimes
\underbrace{
\varpi_1^{u_{j}}\otimes\varpi_1^{v}
}_{}
\otimes\,\varpi_1^{u_{j+1}}\otimes\cdots\otimes\varpi_1^{u_k})
\\
\rightarrow
\widehat{\cal P}&
(\varpi_1^{u_1}\otimes\cdots\otimes\varpi_1^{u_{j-1}}\otimes
\underbrace{\varpi_1^{v}
\otimes\varpi_1^{u_{j}}
}_{}
\otimes\,\varpi_1^{u_{j+1}}\otimes\cdots\otimes\varpi_1^{u_k}
).
\end{align*}
We also put
\begin{align*}
W&\left(\varpi_1^{u_1}\otimes \varpi_1^{u_2}\otimes \cdots\otimes \varpi_1^{u_k},
\varpi_1^{v_1}\otimes \varpi_1^{v_2}\otimes \cdots\otimes \varpi_1^{v_l}\right)\\:=&
\prod_{1\leq j\leq l}^{\longleftarrow}
W\left(\varpi_1^{u_1}\otimes \varpi_1^{u_2}\otimes \cdots\otimes 
\varpi_1^{u_k},\varpi_1^{v_j}\right)^{[j,k+j]} \\ % \label{bfusion}
:\widehat{\cal P}&
(\varpi_1^{u_1}\otimes\cdots\otimes\varpi_1^{u_{k}}\otimes
\varpi_1^{v_1}
\otimes\cdots
\otimes\varpi_1^{v_{l}})
\rightarrow
\widehat{\cal P}
(\varpi_1^{v_1}
\otimes\cdots
\otimes\varpi_1^{v_{l}}
\otimes
\varpi_1^{u_1}\otimes\cdots\otimes\varpi_1^{u_{k}}), 
\end{align*}
where
\begin{align*}
&W(\varpi_1^{u_1}\otimes \cdots\otimes\varpi_1^{u_k},\varpi_1^{v_j})
^{[j,k+j]}
:={\rm id}^{\otimes(j-1)}\otimes
W(\varpi_1^{u_1}\otimes \cdots\otimes\varpi_1^{u_k},\varpi_1^{v_j})
\otimes{\rm id}^{\otimes(l-j)}
\\
:\widehat{\cal P}&
(\varpi_1^{v_1}\otimes\cdots\otimes\varpi_1^{v_{j-1}}\otimes
\underbrace{\varpi_1^{u_1}\otimes\cdots\otimes\varpi_1^{u_k}
\otimes\varpi_1^{v_j}}_{}\otimes\,\varpi_1^{v_{j+1}}\otimes\cdots
\otimes\varpi_1^{v_l}
)\\
\rightarrow
\widehat{\cal P}&
(\varpi_1^{v_1}\otimes\cdots\otimes\varpi_1^{v_{j-1}}\otimes
\underbrace{
\varpi_1^{v_{j}}\otimes
\varpi_1^{u_1}\otimes\cdots\otimes\varpi_1^{u_k}}_{}
\otimes\,\varpi_1^{v_{j+1}}\otimes\cdots\otimes\varpi_1^{v_l}
).
\end{align*}

We will realize the space $\widehat{\cal P}(\varpi_2^u)$ as a subspace
of $\widehat{\cal P}(\varpi_1^{u}\otimes \varpi_1^{u-\h}).$
For this purpose, let us introduce the {\em fusion projector} $\pi_{\varpi_2^u}$
by specializing the parameter in $W\left(\varpi^u_1,\varpi^v_1\right)$:
\begin{equation}
\pi_{\varpi_2^u}:=W (\varpi_1^{u-\h},\varpi_1^u ):
\widehat{\cal P}(\varpi_1^{u-\h}\otimes\varpi_1^u)\rightarrow
\widehat{\cal P}(\varpi_1^u\otimes\varpi_1^{u-\h}).
\label{projector}
\end{equation}
%%%%%%% A basis %%%%%%%%%%
\begin{lem} The space 
$\pi_{\varpi_2^u}(\widehat{\cal P}(\varpi_1^{u-\h}\otimes\varpi_1^u)_\la)$
has a basis $\{\bar{f}_{\la}^{\la+\widehat{r}}(u)\}_{r\in{\cal P}_2}$
given by
\begin{align}
\bar{f}_{\la}^{\la+\widehat{p}+\widehat{q}}(u):= &[\la_{p-q}+\h]
e^{\la+\widehat{p}}_{\la}(u)\otimes 
e^{\la+\widehat{p}+\widehat{q}}_{\la+\widehat{p}}(u-\h)  \nonumber \\
&+[\la_{q-p}+\h]
e^{\la+\widehat{q}}_{\la}(u)\otimes 
e^{\la+\widehat{p}+\widehat{q}}_{\la+\widehat{q}}(u-\h),\label{f1}
\end{align}
where
$p=\pm\ep_1, \, q=\pm\ep_2$ , and 
\begin{equation}
\bar{f}_{\la}^\la(u) := 
\sum_{p \in {\cal P}_1}[2\la_p+2\h] 
e^{\la+\widehat{p}}_{\la}(u)\otimes 
e_{\la+\widehat{p}}^{\la}(u-\h).\label{f2}
\end{equation}
\end{lem}
\begin{prf} For $p,q \in {\cal P}_1$, $q \neq \pm p$, 
we have
\begin{equation*}
\pi_{\varpi_2^u} 
\left(e^{\la+\widehat{p}}_{\la}(u-\h)\otimes 
e^{\la+2\widehat{p}}_{\la+\widehat{p}}(u)\right)
= 
\left( \begin{matrix}  
                &p       &      \\
          p \!\!\!\!    &\boxed{-\h} &\!\!\!\!\! p\\
                &p      &       \\
\end{matrix} \right)
e^{\la+\widehat{p}}_{\la}(u)
\otimes e^{\la+2\widehat{p}}_{\la+\widehat{p}}(u-\h)
=0, 
\end{equation*}
\begin{align*}
&\pi_{\varpi_2^u}
\left(e^{\la+\widehat{p}}_{\la}(u-\h)
\otimes 
e^{\la+\widehat{p}+\widehat{q}}_{\la+\widehat{p}}(u) \right) \\
=& 
\left( \begin{matrix}  
                &p       &      \\
          p \!\!\!\!    &\boxed{-\h} &\!\!\!\!\! q\\
                &q      &       \\
\end{matrix} \right)
e^{\la+\widehat{p}}_{\la}(u)
\otimes e^{\la+\widehat{p}+\widehat{q}}_{\la+\widehat{p}}(u-\h) +
\left( \begin{matrix}  
                &p       &      \\
          q \!\!\!\!    &\boxed{-\h} &\!\!\!\!\! q\\
                &p      &       \\
\end{matrix} \right)
e^{\la+\widehat{q}}_{\la}(u)
\otimes e^{\la+\widehat{p}+\widehat{q}}_{\la+\widehat{q}}(u-\h)  \\
=& \frac{[-2\h]}{[-3\h]\,[\la_{p-q}]} 
\left( [\la_{p-q}+\h] e^{\la+\widehat{p}}_{\la}(u)
\otimes e^{\la+\widehat{p}+\widehat{q}}_{\la+\widehat{p}}(u-\h)
+ [\la_{q-p}+\h] e^{\la+\widehat{q}}_{\la}(u)
\otimes e^{\la+\widehat{p}+\widehat{q}}_{\la+\widehat{q}}(u-\h) \right),
\end{align*}
and
\begin{align*}
&\pi_{\varpi_2^u} \left(e^{\la+\widehat{p}}_{\la}(u-\h)
\otimes e^{\la}_{\la+\widehat{p}}(u) \right) \\
=& 
\sum_{r\in {\cal P}_1}
\left( \begin{matrix}  
                &p       &      \\
          r \!\!\!\!    &\boxed{-\h} &\!\!\!\!\! -p \\
                &-r      &       \\
\end{matrix} \right)
e^{\la+\widehat{r}}_{\la}(u)\otimes e^{\la}_{\la+\widehat{r}}(u-\h)  \\
=&
\frac{[-\h]\,[\la_{p+q}-\h]\,[\la_{p-q}-\h]}{[-3\h]\,[\la_{p+q}]\,[\la_{p-q}]\,[2\la_p]}
\left( \sum_{r \in {\cal P}_1} [2\la_r+2\h] 
e^{\la+\widehat{r}}_{\la}(u)
\otimes e^{\la}_{\la+\widehat{r}}(u-\h)
\right).
\end{align*}
Here we have used the three-term identity (\ref{3term}).
\end{prf}

Thus we know the subspace $\pi_{\varpi_2^u}(\widehat{\cal P}(\varpi_1^{u-\h}\otimes\varpi_1^u)_\la)$
is naturally isomorphic to the space $\widehat{\cal P}(\varpi_2^u)_\la.$
In the following, we will identify the image ${\rm Im}(\pi_{\varpi_2^u})\subset
\widehat{\cal P}(\varpi_1^u\otimes \varpi_1^{u-\h})$ 
with the space $\widehat{\cal P}(\varpi_2^u)$
via
$\bar{f}_\la^\mu(u)\leftrightarrow f_\la^\mu(u).$

\begin{prop}\label{preserve}
Define the operators $\widetilde{W}_{dd'}(u-v)$ by 
\begin{equation}
\widetilde{W}_{21}(u-v) := W(\varpi^u_1 \otimes \varpi^{u-\h}_1, \varpi^v_1),\ 
\widetilde{W}_{12}(u-v) := W(\varpi^u_1 , \varpi^{v}_1 \otimes \varpi^{v-\h}_1)\label{w2112}
\end{equation}
and
\[
\widetilde{W}_{22}(u-v) := W(\varpi^u_1\otimes \varpi^{u-\h}_1 , \varpi^{v}_1 \otimes \varpi^{v-\h}_1).
\]
We have 
\begin{equation*}
 \widetilde{W}_{dd'}(u-v) ( \widehat{\cal P}(\varpi_d^u \otimes \varpi_{d'}^v)_\la^\mu )
\subset \widehat{\cal P}(\varpi_{d'}^v \otimes \varpi_d^u)_\la^\mu. 
\end{equation*}
\end{prop}
\begin{prf}
From the definition of $\pi_{\varpi_2^u}$ (\ref{projector}) and 
the YBE (\ref{YBE1}), 
\begin{align}
&W^{1,2}(u-v) W^{2,3}(u-v-\h) (\pi_{\varpi_2^u}\otimes {\rm id}) \nonumber \\
= &({\rm id} \otimes \pi_{\varpi_2^u}) W^{1,2}(u-v-\h) W^{2,3}(u-v).
\label{W21pres}
\end{align}
Applying this to the definition of $\widetilde{W}_{21}(u-v)$, we get
\begin{equation*}
 \widetilde{W}_{21}(u-v) ( \widehat{\cal P}(\varpi_2^u \otimes \varpi_1^v)_\la^\mu )  
\subset \widehat{\cal P}(\varpi_1^v \otimes \varpi_2^u)_\la^\mu. 
\end{equation*}
By a same argument, we have
\begin{align}
&W^{2,3}(u-v+\h) W^{1,2}(u-v) ({\rm id} \otimes \pi_{\varpi_2^u}) \nonumber \\
= &(\pi_{\varpi_2^u} \otimes {\rm id}) W^{2,3}(u-v) W^{1,2}(u-v+\h),
\label{W12pres}
\end{align}
and
\begin{equation*}
 \widetilde{W}_{12}(u-v) ( \widehat{\cal P}(\varpi_1^u \otimes \varpi_2^v)_\la^\mu )  
\subset \widehat{\cal P}(\varpi_2^v \otimes \varpi_1^u)_\la^\mu. 
\end{equation*}
Together with the equations (\ref{W21pres}),(\ref{W12pres}) and
the definition of $\widetilde{W}_{22}(u-v)$, we obtain
\begin{equation*}
 \widetilde{W}_{22}(u-v) ( \widehat{\cal P}(\varpi_2^u \otimes \varpi_2^v)_\la^\mu )  
\subset \widehat{\cal P}(\varpi_2^v \otimes \varpi_2^u)_\la^\mu. 
\end{equation*}
\end{prf}

We denote by $W_{dd'}(u-v)$ the restricted operators 
$
\widetilde{W}_{dd'}(u-v)|_{\widehat{\cal P}(\varpi_d^u \otimes \varpi_{d'}^v)}
$
and introduce their 
matrix coefficients by the following equation
\begin{eqnarray*}
W_{dd'}(u-v)\,g_\la^\mu(u)\otimes {g}_\mu^\nu(v)
=\sum_{\kappa\in \mathfrak{h}^*} 
W_{dd'}\left(\left.\begin{array}{ll}
          \la   &        \mu        \\
         \kappa &    \nu    \\
\end{array}
\;\right|u-v
\right)\; {g}_\la^\kappa(v)\otimes g_\kappa^\nu(u)
.
\end{eqnarray*}
By the construction,
the operators $W_{dd'}(u-v)$ clearly satisfies
the YBE (\ref{YBE_op}) in operator form, and
their coefficients 
$W_{dd'}\left(\left.\begin{array}{ll}
          \la   &        \mu        \\
         \kappa &    \nu    \\
\end{array}
\;\right|u-v
\right)$
satisfies
the YBE (\ref{YBE1}).
For $p,r \in {\cal P}_d$ and 
$s,q \in {\cal P}_{d'}\;(d,d' = 1,2)$ such that $p+q=r+s$
we write for brevity (as far as confusion does not arise)
\begin{equation}
\begin{matrix}    
                &p              &               \\
   s\!\!\!\!    &\boxed{u}      &\!\!\!\!\! q   \\
                &r              &               \\
\end{matrix}
=
W_{dd'}\left(\left.\begin{array}{ll}
          \la   &        \la+\widehat{p}        \\
         \la+\widehat{s}&\la+\widehat{p}+\widehat{q}    \\
\end{array}
\;\right|u
\right).\label{box_notation}
\end{equation}

We calculate the coefficients of the operator $W_{21}(u)$ as example.
In what follows, we will often omit the dependence of 
$g_\la^\mu(u)\in 
{\widehat{\cal P}}(\varpi_d^u)$
on $u$ (the spectral parameter) for brevity.
Let $p\in {\cal P}_1.$
From the definitions of $f_{\la}^{\la}$ and ${\widetilde W}_{21}$ 
(\ref{f2},\ref{w2112})
we have 
\begin{align*}
&{W}_{21}(u) \,f_{\la}^{\la}\otimes e_{\la}^{\la+\widehat{p}} 
={\widetilde W}_{21}(u) \,f_{\la}^{\la}\otimes e_{\la}^{\la+\widehat{p}} \\
= &{\widetilde W}_{21}(u)  \left( \sum_{r \in {\cal P}_1}
[2\la_{r}+2\h] \,e_{\la}^{\la +\widehat{r}} \otimes e_{\la +\widehat{r}}^{\la} 
\otimes e_{\la}^{\la+\widehat{p}}
\right) \\
=&\sum_{q\in {\cal P}_1}
e_{\la}^{\la+\widehat{q}}\otimes 
\left(
\sum_{\substack{s,t \in {\cal P}_1 \\ s+t=p-q}}
V_q(\la;s,t;u)e_{\la+\widehat{q}}^{\la+\widehat{q}+\widehat{s}} \otimes
e_{\la+\widehat{q}+\widehat{s}}^{\la+\widehat{p}}\right),
\end{align*}
where we denote by $V_q(\la;s,t;u)$ the following function
\begin{equation*}
\sum_{r \in {\cal P}_1}
[2\la_r+2\h] \;W_{11}\left(\left.\begin{array}{ll}
          \la   &        \la+\widehat{r}        \\
         \la+\widehat{q} & \la+\widehat{q}+\widehat{s}    \\
\end{array}
\;\right|u
\right)\,
W_{11}\left(\left.\begin{array}{ll}
          \la+\widehat{r}   &     \la        \\
         \la+\widehat{q}+\widehat{s} & \la+\widehat{p}    \\
\end{array}
\;\right|u-\h
\right) .%\label{exam2}
\end{equation*}

If $q \in {\cal P}_1$ such that $q \neq \pm p$, then the functions
$V_q(\la;s,t;u)$ vanish
except for $(s,t) = (p,-q)\; \mbox{or} \;(-q,p)$,
and one can easily show that 
\begin{equation}
\frac{V_q(\la;p,-q;u)}{[(\la+\widehat{q})_{p+q}+\h]}=
\frac{V_q(\la;-q,p;u)}{[(\la+\widehat{q})_{-q-p}+\h]}.\label{V/[]=V/[]}
\end{equation}
This equation implies that the vector
\begin{equation*}
V_q(\la;p,-q;u)\,e_{\la+\widehat{q}}^{\la+\widehat{q}+\widehat{p}}
\otimes e_{\la+\widehat{q}+\widehat{p}}^{\la+\widehat{p}}+
V_q(\la;-q,p;u) \,e_{\la+\widehat{q}}^{\la}\otimes e_{\la}^{\la+\widehat{p}}
\label{exam}
\end{equation*}
is proportional to 
$f_{\la+\widehat{q}}^{\la+\widehat{p}}$
and its coefficient (the both hands sides of (\ref{V/[]=V/[]}))
is calcurated as
$$
\frac{[u-\h]\,[u+\h]\,[u+3\h]\,[2\h]}{[-3\h]^2\,[\h]^2}
\frac{[\la_{q-p}-\h-u]\,[2\la_q+2\h]}{[\la_{q-p}-\h]\,[\la_{q+p}+\h]}
$$
by using 
the three term identity (\ref{3term}). 
This function is labeled by (see (\ref{box_notation}))
$$\begin{matrix}  
                &0       &      \\
          q \!\!\!\!    &\boxed{u} &\!\!\!\!\! p\\
                &p-q      &       \\
\end{matrix} \quad(q\ne\pm p).
$$

Let us consider the term for $q=p.$ 
For all $s\in {\cal P}_1$ we have from the three term identity
\begin{equation}
\frac{V_p(\la;s,-s;u)}{[2(\la+\widehat{p})_s+2\h]}
=\frac{[u-\h]\,[u+\h]\,[u+3\h]}{[-3\h]^2 \,[\h]}
\frac{[u+\h]}{[\h]}
\prod_{\substack{r \in {\cal P}_1\\ r\neq \pm p}} 
\frac{[\la_{p+r}+2\h]}{[\la_{p+r}+\h]}.\label{ex1'}
\end{equation}
The right hand side of this equation is independent of 
$s\in{\cal P}_1.$
Thus we see that the vector 
$$
\sum_{s\in{\cal P}_1}
V_p(\la;s,-s;u)e_{\la+\widehat{p}}^{\la+\widehat{p}+\widehat{s}} \otimes
e_{\la+\widehat{p}+\widehat{s}}^{\la+\widehat{p}}
$$
is proportional to $f_{\la+\widehat{p}}^{\la+\widehat{p}}$
and its coefficient is equal to the right hand side of (\ref{ex1'}),
which is labeled by 
$$\begin{matrix}  
                &0       &      \\
          p \!\!\!\!    &\boxed{u} &\!\!\!\!\! p\\
                &0      &       \\
\end{matrix}.
$$

Here we write all fused Boltzmann weights 
(the coefficients of the operator $W_{21}(u)$).
They are obtained by the three term identity (\ref{3term}).
We assume $p,q \in {\cal P}_1$ satisfy $p \neq \pm q$.
The common factor $[u-\h]\,[u+\h]\,[u+3\h][-3\h]^{-2}[\h]^{-1}$ is dropped.
\begin{align}
\begin{matrix}  
                &p+q       &      \\
          q \!\!\!\!    &\boxed{u} &\!\!\!\!\! q\\
                &p+q     &       \\
\end{matrix}
&=
\frac{[u+2\h]}{[\h]},\nonumber\\
\begin{matrix}  
                &p-q       &      \\
          q \!\!\!\!    &\boxed{u} &\!\!\!\!\! q\\
                &p-q     &       \\
\end{matrix}
&=
\frac{[u]}{[\h]}
\frac{[2\la_q+2\h]}{[2\la_q]}
\frac{[\la_{p-q}-\h]}{[\la_{p-q}+\h]},\nonumber\\
\begin{matrix}  
                &0       &      \\
          q \!\!\!\!    &\boxed{u} &\!\!\!\!\! q\\
                &0     &       \\
\end{matrix}
&=
\frac{[u+\h]}{[\h]}
\prod_{\substack{r \in {\cal P}_1\\ r\neq \pm q}} \frac{[\la_{q+r}+2\h]}{[\la_{q+r}+\h]},\label{ex1} \\
\begin{matrix}  
                &q-p       &      \\
          q \!\!\!\!    &\boxed{u} &\!\!\!\!\! p\\
                &0      &       \\
\end{matrix}
&=
\frac{[\la_{q-p}-u]\,[\la_{q+p}+2\h]}{[2\la_p]\,[\la_{q-p}+\h]}, \label{ex2} \\
\begin{matrix}  
                &0       &      \\
          q \!\!\!\!    &\boxed{u} &\!\!\!\!\! p\\
                &p-q      &       \\
\end{matrix}
&=
\frac{[2\h]}{[\h]}
\frac{[\la_{q-p}-\h-u]\,[2\la_q+2\h]}{[\la_{q-p}-\h]\,[\la_{q+p}+\h]}, \label{ex3} \\
\begin{matrix}  
                &p+q       &      \\
          q \!\!\!\!    &\boxed{u} &\!\!\!\!\! -q\\
                &p-q      &       \\
\end{matrix}
&=
\frac{[2\h]}{[\h]}
\frac{[2\la_q-u]\,[\la_{p-q}-\h]}{[2\la_q]\,[\la_{p+q}+\h]}.\nonumber
\end{align}

Next we give the example of $W_{12}$. 
In this case, the common factor $[u]\,[u+2\h]\,[u+4\h]\,[-3\h]^{-2}[\h]^{-1}$ 
is dropped.
To obtain them, we use only the three-term identity (\ref{3term}).
\begin{align}
\begin{matrix}  
                &p       &      \\
          p\! +\! q \!\!\!\!    &\boxed{u} &\!\!\!\!\! p \! + \! q\\
                &p     &       \\
\end{matrix}
&=
\frac{[u+3\h]}{[\h]},\nonumber\\
\begin{matrix}  
                      &p       &      \\
          q\! -\!p \!\!\!\!    &\boxed{u} &\!\!\!\!\! q\! - \! p\\
                &p     &       \\
\end{matrix}
&=
\frac{[u+\h]}{[\h]}
\frac{[2\la_p-2\h]}{[2\la_p]}
\frac{[\la_{q-p}+2\h]}{[\la_{q-p}]},\nonumber\\
\begin{matrix}  
                &p       &      \\
          0 \!\!\!\!    &\boxed{u} &\!\!\!\!\! 0\\
                &p     &       \\
\end{matrix}
&=
\frac{[u+2\h]}{[\h]}
\prod_{\substack{r\in {\cal P}_1 \\ r\neq \pm p}}
\frac{[\la_{p+r}-\h]}{[\la_{p+r}]},\label{sign1}\\
\begin{matrix}  
                &p       &      \\
          0 \!\!\!\!    &\boxed{u} &\!\!\!\!\! q \! - \! p\\
                &q      &       \\
\end{matrix}
&=
\frac{[\la_{p-q}-2\h-u]\,[\la_{p+q}-\h]}{[2\la_p]\,[\la_{q-p}]},\nonumber\\
\begin{matrix}  
                &p       &      \\
          p \! - \! q \!\!\!\!    &\boxed{u} &\!\!\!\!\! 0\\
                &q      &       \\
\end{matrix}
&=
\frac{[2\h]}{[\h]}
\frac{[\la_{p-q}-\h-u]\,[2\la_q-2\h]}{[\la_{q-p}]\,[\la_{p+q}]},\nonumber\\%\label{sign2}
\begin{matrix}  
                &p       &      \\
          p \! + \! q \!\!\!\!    &\boxed{u} &\!\!\!\!\! q\! -\! p \\
                &-p      &       \\
\end{matrix}
&=
\frac{[2\h]}{[\h]}
\frac{[2\la_p-\h-u]\,[\la_{p+q}+2\h]}{[2\la_p]\,[\la_{q-p}]}.\nonumber\end{align}

Finally we give the example of $W_{22}$. 
They are equivalent to the Boltzmann weights associated 
to the vector representation of the type $B_2$ Lie algebra
(see Ref.\cite{JMO2}).
We write only two cases as example, which is used to define the
difference operator $M_2(u)$.        
We will drop the common factor $G(u)$ (\ref{G(u)}) here.
\begin{align}
\begin{matrix}  
                &p+q       &      \\
          0 \!\!\!\!    &\boxed{u} &\!\!\!\!\! 0\\
                &p+q      &       \\
\end{matrix}
&=
\frac{[\la_{p+q}-\h]}{[\la_{p+q}+\h]},\label{deg2part}
\\
\begin{matrix}  
                &0      &      \\
          0 \!\!\!\!    &\boxed{u} &\!\!\!\!\! 0\\
                &0      &       \\
\end{matrix}
&=
\frac{[2\h]}{[6\h]}
\left(
\sum_{\substack{r=\pm\ep_1 \\ s=\pm\ep_2}}
\frac{[2\la_r+2\h][2\la_s+2\h]}{[2\la_r][2\la_s]}
\frac{[\la_{r+s}-5\h][\la_{r+s}+2\h]}{[\la_{r+s}][\la_{r+s}+\h]}
-\frac{[u+6\h]\,[u-3\h]}{[u]\,[u+3\h]}
\right).\label{degzero}
\end{align}
The formulas (\ref{sign1}), (\ref{deg2part}) and (\ref{degzero})
together give the explicit form of 
$\widetilde{M}_d$ (Theorem \ref{commute} (ii)).

We explain how to calculate the fused Boltzmann weight 
$
\begin{matrix}  
                &0      &      \\
          0 \!\!\!\!    &\boxed{u} &\!\!\!\!\! 0\\
                &0      &       \\
\end{matrix}
$.
According to the definition of the operator $W_{22}(u)$ and the vector 
$f_\la^\la$ (\ref{f2}), the 
coefficient of $W_{22}(u) f_\la^\la \otimes f_\la^\la$ with respect to
$f_\la^\la \otimes f_\la^\la$ is
equal to
\begin{equation}
\frac{1}{[2\la_p+2\h]}
\sum_{r \in {\cal P}_1}
[2\la_r+2\h]\; W_{21}\left(\left.\begin{array}{ll}
          \la   &        \la        \\
         \la+\widehat{p} & \la+\widehat{r}    \\
\end{array}
\;\right|u
\right)
W_{21}\left(\left.\begin{array}{ll}
          \la+\widehat{p}   &  \la+\widehat{r}     \\
          \la           &  \la    \\
\end{array}
\;\right|u+\h
\right) .
\label{00001}
\end{equation}
In this summation, if $r$ is equal to $-p$, then
\[
W_{21}\left(\left.\begin{array}{ll}
          \la   &        \la        \\
         \la+\widehat{p} & \la-\widehat{p}    \\
\end{array}
\;\right|u
\right) = 0. 
\]
So that (\ref{00001}) can be rewritten as
\begin{align*}
&W_{21}\left(\left.\begin{array}{ll}
          \la   &        \la        \\
         \la+\widehat{p} & \la+\widehat{p}    \\
\end{array}
\;\right|u
\right)
W_{21}\left(\left.\begin{array}{ll}
          \la+\widehat{p}   &  \la+\widehat{p}     \\
          \la           &  \la    \\
\end{array}
\;\right|u+\h
\right)   \\
{}+
\sum_{\substack{q \in {\cal P}_1 \\ q \neq \pm p}}  
&\frac{[2\la_q+2\h]}{[2\la_p+2\h]}\;
W_{21}\left(\left.\begin{array}{ll}
          \la   &        \la        \\
         \la+\widehat{p} & \la+\widehat{q}   \\
\end{array}
\;\right|u
\right)
W_{21}\left(\left.\begin{array}{ll}
          \la+\widehat{p}   &   \la+\widehat{q}        \\
          \la           &   \la    \\
\end{array}
\;\right|u+\h
\right).
\end{align*}
By means of (\ref{ex1}), (\ref{ex2}) and (\ref{ex3}),
this function is equal to
\begin{align*}
&\frac{[u-\h]\,[u]\,[u+\h]\,[u+2\h]\,[u+3\h]\,[u+4\h]\,[2\h]}{[-3\h]^3\,[\h]^4}\left(
\frac{[u+\h]\,[u+2\h]}{[2\h]\,[-3\h]}
\prod_{\substack{q \in {\cal P}_1 \\ q \neq \pm p}}
\frac{[\la_{p+q}-\h][\la_{p+q}+2\h]}{[\la_{p+q}][\la_{p+q}+\h]}\right. \\
&+ \left. \frac{[\h]}{[-3\h]}\sum_{\substack{q\in {\cal P}_1 \\ q \neq \pm p}}
\frac{[2\la_q-2\h]}{[2\la_q]}
\frac{[\la_{p+q}+2\h+u]\,[\la_{p+q}-\h-u]\,[\la_{p-q}-\h]}
{[\la_{p+q}]\,[\la_{p+q}-\h]\,[\la_{p-q}+\h]} \right) .
\end{align*}
To obtain the formula (\ref{degzero}), we use the following lemma.
\begin{lem} For any $p\in {\cal P}_1$, we have
\begin{align}
&\frac{[u+\h]\,[u+2\h]}{[2\h]\,[-3\h]}
\prod_{\substack{q\in {\cal P}_1 \\ q \neq \pm p}}
\frac{[\la_{p+q}-\h]}{[\la_{p+q}]}
\frac{[\la_{p+q}+2\h]}{[\la_{p+q}+\h]} \nonumber\\
&{}+\frac{[\h]}{[-3\h]}
\sum_{\substack{q\in {\cal P}_1 \\ q \neq \pm p}}
\frac{[2\la_q-2\h]}{[2\la_q]}
\frac{[\la_{p+q}+2\h+u]\,[\la_{p+q}-\h-u]\,[\la_{p-q}-\h]}
{[\la_{p+q}]\,[\la_{p+q}-\h]\,[\la_{p-q}+\h]}\nonumber\\
=&\frac{[u][u+3\h]}{[6\h][-3\h]}
\sum_{\substack{r=\pm\ep_1 \\ s=\pm\ep_2}}
\frac{[2\la_r+2\h][2\la_s+2\h]}{[2\la_r][2\la_s]}
\frac{[\la_{r+s}-5\h][\la_{r+s}+2\h]}{[\la_{r+s}][\la_{r+s}+\h]}
+\frac{[u+6\h][u-3\h]}{[6\h][u+3\h]}.
\label{eq:ZeroDegTerm}
\end{align}
\end{lem}
\begin{prf}
Let $f(\la_p)$ be (the left-hand side) $-$ (the right-hand side)
of (\ref{eq:ZeroDegTerm}), regarded as a function of $\la_p$.
It is doubly periodic function of the periods $1,\tau$.
Let us show that it is entire.
The apparent poles of $f(\la_p)$ are located at
$$\la_p=\la_q,\;\la_p=\la_q\pm \h\,
(p,q\in {\cal P}_1,\;p+q\ne 0),\;\la_p=0\,(p\in {\cal P}_1).$$
Note that the left-hand side of (\ref{eq:ZeroDegTerm}) 
is clearly invariant under 
$\la_q\mapsto-\la_q$, and the right-hand side is $W$-invariant.
In view of the symmetry, it suffices to check the regularity
at $\la_p=\la_q,\;\la_p=\la_q- \h$ and $\la_p=0.$ 
By the three-term identity (\ref{3term}), it is easy to see that 
the residue of $f(\la_p)$ at $\la_p=\la_q- \h$ vanishes.
Manifestly, the point $\la_p=\la_q$ and $\la_p=0$ is regular.

Now we have proved that $f(\la_p)$ is independent of $\la_p$. 
We will show $f(-\la_q-2\h)=0$.
This can be directly checked by using
the identity (\ref{3term}) twice,
and the proof completes.
\end{prf}

\section{Commutativity of the difference operators}
\label{difference}
This section is devoted to the proof of commutativity of 
the difference operators(Theorem\ref{commute} (i)).
For $t\in {\cal P}_d+{\cal P}_{d'}$ we will introduce the matrices 
$A_t(\la|u,v),B_t(\la|v,u)$
whose index set is $I_{t}:=\{(p,q)\in{\cal P}_d\times{\cal P}_{d'}\;|\;p+q=t\}:$

\begin{equation*}
A_t(\la|u,v)^{(p,q)}_{(r,s)}:=
W_{d2}\left(\left.\begin{array}{ll}
          \la   &        \la+\widehat{p}        \\
          \la   &        \la+\widehat{r}   \\
\end{array}
\;\right|u
\right)
W_{d'2}\left(\left.\begin{array}{ll}
          \la+\widehat{p}   &        \la+\widehat{t}        \\
          \la+\widehat{r}   &        \la+\widehat{t}   \\
\end{array}
\;\right|v
\right),
\end{equation*}
\begin{equation*}
B_t(\la|v,u)^{(p,q)}_{(r,s)}
:=
W_{d'2}\left(\left.\begin{array}{ll}
          \la   &        \la+\widehat{q}        \\
          \la   &        \la+\widehat{s}   \\
\end{array}
\;\right|v
\right)
W_{d2}\left(\left.\begin{array}{ll}
          \la+\widehat{q}   &        \la+\widehat{t}        \\
          \la+\widehat{s}   &        \la+\widehat{t}   \\
\end{array}
\;\right|u
\right).
\end{equation*}
With these matrices, we can write down both the left and right hand sides as
$$M_d(u)M_{d'}(v)=\sum_{t\in {\cal P}_d+{\cal P}_{d'}}{\rm tr}\,A_t(\la|u,v)\;T_{\widehat{t}},
\quad M_{d'}(v)M_d(u)=\sum_{t\in {\cal P}_d+{\cal P}_{d'}}{\rm tr}\,B_t(\la|v,u)\;T_{\widehat{t}}.$$
Let us also define the matrix $W_t(\la|u-v)$ with the same index set:
\begin{equation*}
W_t(\la|u-v)^{(p,q)}_{(r,s)}
:=
W_{dd'}\left(\left.\begin{array}{ll}
          \la   &        \la+\widehat{p}        \\
          \la+\widehat{s}    &        \la+\widehat{t}   \\
\end{array}
\;\right|u-v
\right).
\end{equation*}

The YBE (\ref{YBE1}) implies
\begin{equation*}
W_t(\la|u-v) A_t(\la|u,v) 
=B_t(\la|v,u) W_t(\la|u-v).
\end{equation*}
By the inversion relation (\ref{FirstInversion}),
it can be seen that $W_t(\la|u-v)$ is invertible for 
generic $u,v\in {\mathbb C}.$
It follows that 
${\rm tr}\,A_t(\la|u,v)={\rm tr}\,B_t(\la|v,u)$
for all $u,v \in {\mathbb C}$.
Hence we have $M_d(u)M_{d'}(v)=M_{d'}(v)M_d(u)$ for all 
$u,v \in {\mathbb C}.$

\section{Space of Weyl group invariant theta functions}\label{Theta}
This section is devoted to the proof of Theorem \ref{theta}.
Let $Q^\vee,P^\vee$ be the coroot and coweight lattice respectively.
Under the identification $\mathfrak{h}=\mathfrak{h}^*$ via the from $(\,,\,),$
these are given by 
$$Q^\vee=\itg 2\ep_1\oplus\itg 2\ep_2,
\quad P^\vee=Q^\vee
+\itg(\ep_1+\ep_2).
$$
\begin{lem}For all $\be\in P^\vee$ and $d=1,2$, we have 
\begin{equation}
[S_{\tau\beta},M_d(u)]=[S_{\beta},M_d(u)]=0.
\label{SM=MS}
\end{equation}
\end{lem}
\begin{prf}
Note that if $p,q\in {\cal P}_1\;(q\ne\pm p)$ 
and $\beta\in P^{\vee}$ then $\beta_{p+q} \in \itg$.
By the 
quasi-periodicity (\ref{quasi_Per}), we have
\begin{equation*}
\frac{[(\la+\tau\beta)_{p+q}-\h]}
{[(\la+\tau\beta)_{p+q}]}=
e^{2\pi i\beta_{p+q}\h}
\frac{[\la_{p+q}-\h]}
{[\la_{p+q}]},\quad
\frac{[(\la+\beta)_{p+q}-\h]}
{[(\la+\beta)_{p+q}]}=
\frac{[\la_{p+q}-\h]}
{[\la_{p+q}]}.
\end{equation*}
Using these equations, we have for all $p\in{\cal P}_1$
\begin{align*}
S_{\tau\beta}\prod_{q\ne \pm p}
\frac{[\la_{p+q}-\h]}{[\la_{p+q}]}
T_{\widehat{p}}f(\la)
&=
e^{2\pi i((\la,\beta)+\tau(\beta,\beta)/2)}
\prod_{q\ne \pm p}
\frac{[(\la+\tau\beta)_{p+q}-\h]}{[(\la+\tau\beta)_{p+q}]}
f(\la+\tau\beta+{\widehat{p}}) \\
&=
e^{2\pi i((\la,\beta)+\tau(\beta,\beta)/2+2\beta_p \h)}
\prod_{q\ne \pm p}
\frac{[\la_{p+q}-\h]}{[\la_{p+q}]}
f(\la+\tau\beta+{\widehat{p}})  \\
&=\prod_{q\ne \pm p}
\frac{[\la_{p+q}-\h]}{[\la_{p+q}]}
T_{\widehat{p}}S_{\tau\beta}f(\la),
\end{align*}
and 
\begin{align*}
&S_{\beta}\prod_{q\ne \pm p}
\frac{[\la_{p+q}-\h]}{[\la_{p+q}]}
T_{\widehat{p}}f(\la)
=
\prod_{q\ne \pm p}
\frac{[(\la+\beta)_{p+q}-\h]}{[(\la+\beta)_{p+q}]}
f(\la+\beta+{\widehat{p}})\nonumber\\
=&
\prod_{q\ne \pm p}
\frac{[\la_{p+q}-\h]}{[\la_{p+q}]}
f(\la+\beta+{\widehat{p}})
=\prod_{q\ne \pm p}
\frac{[\la_{p+q}-\h]}{[\la_{p+q}]}
T_{\widehat{p}}S_{\beta}f(\la).
\end{align*}
Note that $2\beta_p \h=(\widehat{p},\beta)$ etc.
Hence we have $[S_{\tau\beta},M_1(u)]=[S_{\beta},M_1(u)]=0.$
In the same way, we can see that 
the principal part
of ${\widetilde M}_2$
commutes with 
$S_{\tau\beta}$ and $S_{\beta},$ using the equations
\begin{equation*}
\frac{[(\la+\tau\beta)_{p+q}-\h]}
{[(\la+\tau\beta)_{p+q}+\h]}=
e^{2\pi i(2\beta_{p+q}\h)}
\frac{[\la_{p+q}-\h]}
{[\la_{p+q}+\h]},\quad
\frac{[(\la+\beta)_{p+q}-\h]}
{[(\la+\beta)_{p+q}+\h]}=
\frac{[\la_{p+q}-\h]}
{[\la_{p+q}+\h]}.
\end{equation*}
Using (\ref{quasi_Per}) it is easy to see 
that the function
$$
C_{p,q}(\la):=\frac{[2\h]}{[6\h]}
\frac{[2\la_p+2\h]}{[2\la_p]}
\frac{[2\la_q+2\h]}{[2\la_q]}
\frac{[\la_{p+q}-5\h]}{[\la_{p+q}+\h]}
\frac{[\la_{p+q}+2\h]}{[\la_{p+q}]}\;(p,q\in {\cal{P}}_1,p+q\ne 0)
$$ satisfies 
$C_{p,q}(\la+\beta)=C_{p,q}(\la+\tau\beta)
=C_{p,q}(\la)\;(\forall \beta\in P^\vee).$
This means that $S_{\tau \be},\;S_{\be}(\beta\in P^\vee)$ commute with
a multiplication by $C_{p,q}(\la).$ 
\end{prf}

\begin{lem}
For all $\gamma\in P^{\vee}$, we have
\begin{equation}
S_{\tau\gamma}Th^W\subset Th^W,
S_{\gamma}Th^W\subset Th^W.\label{SPreserveTh}
\end{equation}
\end{lem}
\begin{prf}
Let $f\in Th^W$ and $\gamma\in P^{\vee}.$
Since the bilinear form $(\,,\,)$ is $W$-invariant, we have
$(S_{\tau \gamma}f)(w\lambda)=(S_{\tau w^{-1}(\gamma)}f)(\lambda).$
Using (\ref{Heise}), we can write this as $(S_{\tau\gamma}S_{\tau(w^{-1}(\gamma)-\gamma)}f)(\lambda),$
which is equal to 
$S_{\tau\gamma}f(\lambda)$
in view of $w^{-1}(\gamma)-\gamma\in Q^\vee.$
In the same way, we can show that $(S_{\gamma}f)(w\lambda)=(S_{\gamma}f)(\lambda).$

Evidently $S_{\tau\gamma}f$ and $S_{\gamma}f$ are
holomorphic.
For all $\alpha\in Q^\vee,$
using (\ref{Heise}) and $(\gamma,\alpha)\in\mathbb Z,$
it can be seen that
the operators
$S_{\alpha}, S_{\tau\alpha}$
commute with $S_{\gamma},S_{\tau\gamma}.$
Hence $S_{\tau\gamma}f$ or $S_{\gamma}f$ are fixed 
by $S_{\tau\alpha}$ and $S_{\alpha}.$
\end{prf}

Here we prove Theorem \ref{theta}.\

\begin{prf2}
Let $f$ be any function in $Th^W.$
In view of (\ref{SM=MS}), we have 
$S_\al\widetilde{M}_df=S_{\tau\al}\widetilde{M}_df=\widetilde{M}_df$
for all $\al\in Q^\vee\subset P^\vee$.
It is clear from the explicit form of $\widetilde{M}_d$
that $\widetilde{M}_df(w\la)=\widetilde{M}_df(\la)$
for all $w\in W$.

Let us show that the function $\widetilde{M}_df$ is holomorphic
on $\mathfrak{h}^*$. 
For $\mu\in \mathfrak{h}^*$ and $z\in \cpx$, 
we denote by $D_{\mu}^z$ the line in $\mathfrak{h}^*$ defined by
$$D_{\mu}^z:=\{\la\in\mathfrak{h}^*\,|\,(\la,\mu)+z=0\}.$$
The coefficients of the difference operators $\widetilde{M}_d$
have their possible simple poles along
$D+P^\vee+\tau P^\vee$,
where we put 
$$D:=\bigcup_{p\in R_+}D_{p}^0\cup
\bigcup_{q\in{\cal P}_2-\{0\}}D_{q}^\h$$
and $R_+$ is a fixed set of positive roots.

Next we will show that for any
function $f$ in $Th^W$, 
$\widetilde{M}_df$ is regular along $D.$
Let us consider the meromorphic function
$g := \left( \prod_{p\in R_+} [\la_p] \right) \,\widetilde{M}_df$,
which is regular along $D^0:=\bigcup_{p\in R_+}D_{p}^0.$
Since $\widetilde{M}_df$ is $W$-invariant,
it is clear that $g$ is $W$-anti-invariant.
This implies that $g$ has zero along $D^0$
and hence $\widetilde{M}_df$ is regular along $D^0.$

The holomorphy along $\bigcup_{q\in{\cal P}_2-\{0\}}D_{q}^\h$
is somewhat nontrivial. 
Let $p=\pm\ep_1,\,q=\pm\ep_2$.
Clearly, $\widetilde{M}_1f$ is regular along $D_{p+q}^\h.$
Let us consider the function $\widetilde{M}_2f.$
It suffices to show that the following function is regular
along $D_{p+q}^\h$:
$$\frac{[\la_{p+q}-\h]}{[\la_{p+q}+\h]}
T_{\widehat{p}}T_{\widehat{q}}f(\la)
+\frac{[2\h]}{[6\h]}
\frac{[2\la_p+2\h]}{[2\la_p]}
\frac{[2\la_q+2\h]}{[2\la_q]}
\frac{[\la_{p+q}-5\h]}{[\la_{p+q}+\h]}
\frac{[\la_{p+q}+2\h]}{[\la_{p+q}]}
f(\la).$$
We note that, for any $W$-invariant
function $f$, we have 
$\left(T_{\widehat{p}}T_{\widehat{q}}f-f\right)|_{D_{p+q}^\h}=0.$
In view of this, the residue of the above function along ${D_{p+q}^\h}$ is 
easily seen to vanish.
Thus we have proved that for any
function $f$ in $Th^W,$ the functions 
$\widetilde{M}_df(d=1,2)$ are regular along $D.$

For $\beta,\gamma\in P^\vee$, we have, by the definitions 
of $S_{\tau\beta},S_{\gamma}$ and (\ref{SM=MS}), 
\begin{align}
\widetilde{M}_df(\la+\be\tau+\gamma)
&=e^{-\PI((\la,\be)+\tau(\be,\be)/2)}
S_{\tau\be}S_\gamma \widetilde{M}_df(\la) \nonumber \\
&=e^{-\PI((\la,\be)+\tau(\be,\be)/2)}
\widetilde{M}_dS_{\tau\be}S_\gamma f(\la).
\label{ArgumentOfTranslation}
\end{align}
Since $S_{\tau\be}S_\gamma f$ 
belongs to $Th^W$ by (\ref{SPreserveTh}),
$\widetilde{M}_dS_{\tau\be}S_\gamma f$ 
is regular along $D.$
Then (\ref{ArgumentOfTranslation}) implies that
$\widetilde{M}_df$ is regular along $D+\be\tau+\gamma.$
The proof is completed.
\end{prf2}

\section*{Acknowledgments}
The authors are grateful to 
Masato Okado and Michio Jimbo who kindly let them know
detailed points in the proof of the YBE of the 
Boltzmann weight. They also thank Toshiki Nakashima for
useful information and Gen Kuroki,
Masatoshi Noumi, Hiroyuki Ochiai, Yasuhiko Yamada , 
Yasushi Komori and Kazuhiro Hikami
for fruitful discussions and kind interest.
Much of the work was done when T.I. was a postdoctral student at
Mathematical Institute, Tohoku University. He would
like to thank the staff of the Institute for their
support.
%%%%%%%%%%%%% Appendix %%%%%%%%%%%%

%%%%%%%%%%%%% quick hack for JMP appendix style
\section*{Appendix: Similarity transformation}
\def\thesection{A}
\setcounter{equation}{0}
%%%%%%%%%%%%%

Our Boltzmann weights in section \ref{DefOfFace} and the original form in Ref.\cite{JMO2} 
are slightly different.
The original form of type (\ref{2-5c}) and (\ref{2-5d}) are given as follows:
\begin{align}
\begin{matrix}    
                &q       &      \\
          p \!\!\!\!     &\boxed{u} &\!\!\!\!\! p\\
                &q      &       \\
\end{matrix}&=
\frac{[c-u]\,[u]}{[c]\,[\h]}\left( \frac{[\la_{p-q}+\h]\,[\la_{p-q}-\h]}{[\la_{p-q}]^2}\right)^{1/2}
\qquad &(p\ne \pm q),  \label{2-5oc}
\\
\begin{matrix}    
                &q       &      \\
          p \!\!\!\!    &\boxed{u} &\!\!\!\!\! -q\\
                &-p &       \\
\end{matrix}&=
\frac{[u]\,[\la_{p+q}+\h+c-u]}{[c]\,[\la_{p+q}+\h]}
\left( G_{\la p} G_{\la q} \right)^{1/2}
\qquad &(p \ne q). \label{2-5od}
\end{align}
All the other Boltzmann weights ((\ref{2-5a}),(\ref{2-5b}),(\ref{2-5e})) 
are the same as the ones 
we adopted in the section \ref{DefOfFace}.
We denote these weights by 
$W_{JMO}\!
\left(\left.\begin{array}{ll}
        \la     &       \mu     \\
        \kappa  &       \nu     \\
\end{array}\,\right|u
\right) $.
Our Boltzmann weights are obtained from those by the following way.
We introduce an ordering
on the set ${\cal P}$
as
$$\ep_1\prec\ep_2\prec\dots\prec\ep_n
\prec-\ep_n\prec\dots\prec-\ep_2\prec-\ep_1.
$$
For $\la,\mu \in \mathfrak{h}^*$, such that
$\mu - \la = \widehat{q}\in 2\h{\cal P}$,
we define the function $s(\la,\mu)$ by
\begin{equation}
s(\la,\mu)
:=\prod_{\substack{p \in {\cal P} \\ p\prec q}}
[\la_{p-q}]^{-1/2}[\mu_{p-q}]^{-1/2}.
\label{function_s}
\end{equation}
The relation between the Boltzmann weights $W$ in the section \ref{DefOfFace} and the ones in Ref.\cite{JMO2}
is as follows: 
\begin{equation}W\!
\left(\left.\begin{array}{ll}
        \la     &       \mu     \\
        \kappa  &       \nu     \\
\end{array}\,\right|u
\right)
=\frac{s(\la,\mu)s(\mu,\nu)}{s(\la,\kappa)s(\kappa,\nu)}
W_{JMO}\!
\left(\left.\begin{array}{ll}
        \la     &       \mu     \\
        \kappa  &       \nu     \\
\end{array}\,\right|u
\right).\label{Rel_to_JMO}
\end{equation}

\end{document}